\documentclass{tbimc}
\usepackage{graphicx}
\newcommand{\goto}{\rightarrow}

\newcommand{\bff}[1]{{\mbox{\boldmath$#1$}}}
\def\for{\forall} \def\mc{\mathcal} 
\def\mS{{\mathcal S}}

\def\mS{{\mathcal S}}

\def\mT{{\mathcal T}}  
  
 \def\mCB F{${\mathcal CB F}$}

\def\Rui{\Psi}  
 \def\Rui{\Psi}

\def\mG{\mathcal G} 
\newcommand{\be}{\begin{eqnarray}}
\newcommand{\ee}{\end{eqnarray}}
\newcommand{\ba}{\begin{array}}
\newcommand{\ea}{\end{array}}
\newcommand{\baa}{\left[\begin{array}}
\newcommand{\eaa}{\end{array}\right]}

\def\BEN{\begin{enumerate}}  \def\BI{\begin{itemize}}
\def\EEN{\end{enumerate}}   \def\EI{\end{itemize}}

\def\beq{\begin{eqnarray}} \def\eeq{\end{eqnarray}}
\def\bea{\begin{eqnarray*}}
\def\eea{\end{eqnarray*}}
\def\le{\left} \def\ri{\right}

\def\T{\widetilde}  \def\H{\widehat}
\def\bar{\overline}

\def\I{\infty}  \def\a{\alpha}

\def\g{\gamma}  \def\d{\delta} \def\de{\delta}
  \def\th{\theta}
\def\e{\epsilon} \def\k{\kappa} \def\l{\lambda}  
\def\lm{\Pi}

\def\x{\xi}
  \def\r{\rho} \def\s{\sigma}
\def\t{\tau}      
   \def\q{q} \def\qu{\quad} 
\def\F{\Phi}   \def\O{\Omega} 
  
\def\td{\text{\rm d}}

\def\la{\label} \def\fr{\frac} \def\im{\item}

\def\mbb{\mathbb} \def\mbf{\mathbf} 
\def\mc{\mathcal}   \def\bsq{$\blacksquare$}

\newtheorem{Thm}{Theorem}

\def\beT{\begin{Thm}
  }
  \def\eeT{\end{Thm}}

  \newtheorem{Qu}{Problem}
\def\beQ{\begin{Qu}}
  \def\eeQ{\end{Qu}}

\newtheorem{Lem}{Lemma}
\newtheorem{Exa}{Example}

\newtheorem{Cor}{Corollary}

\newtheorem{Def}{Definition}

\newtheorem{Rem}{Remark}

\newtheorem{Exe}{Exercice}
\def\beXe{\begin{Exe}} \def\eeXe{\end{Exe}}

\def\eeD{\end{Def}} \def\beD{\begin{Def}}
\def\beXa{\begin{Exa}} \def\eeXa{\end{Exa}}
\def\beR{\begin{Rem}} \def\eeR{\end{Rem}}
\def\beL{\begin{Lem}} \def\eeL{\end{Lem}}
\newtheorem{Pro}{Proposition}
\def\beP{\begin{Pro}} \def\eeP{\end{Pro}}
\def\beC{\begin{Cor}}
\def\eeC{\end{Cor}}
\def\bc{\begin{cases}
  }
\def\ec{\end{cases}}
\def\BEN{\begin{enumerate}} \def\EEN{\end{enumerate}}

\long\def\symbolfootnote[#1]#2{
\begingroup
\def\thefootnote{\fnsymbol{footnote}}\footnote[#1]{#2}
\endgroup}
\def\fn{\symbolfootnote}
\def\ssec{\subsection} \def\sec{\section}

\def\no{\nonumber} 
\def\bep{\begin{pmatrix}}
  \def\eep{\end{pmatrix}}
  \def\Eq{\Leftrightarrow}
  \def\CL{Cram\'er-Lundberg }

  \def\und{\underline}
  
    \def\P{{\mathbb P}}
    \def\E{{\mathbb E}}
     \def\R{{\mathbb R}}
      
    \newcommand{\1}{\mathop{\mathrm{I}}}

 \def\Lra{\Longrightarrow}

\def\nne{nonnegative }

\def\1{\mathbbm{1}}

\def\1{\mathbf{1}}
\def\mbb{\mathbb}
\def\mI{\mathcal I}

\def\sn{spectrally negative }
\def\sp{spectrally positive }
 \def\LT{Laplace transform}

\def\kil{\mathcal E}

\newcommand{\red}{\textcolor[rgb]{1.00,0.00,0.00}}

\def\qu{\quad}    
\def\lev{L\'evy }   \def\rp{ruin probability }  
\def\deF{de Finetti }  
\def\GS{Gerber-Shiu }
\newcommand{\diff}{{\rm d}} \def\fp{first passage }
\def\ovl{\overline}
    \def\Eb{\E^{b]}}   \def\Ez{\E^{[0}}
    \def\Ezb{\E^{[0,b]}}
    \def\Eazb{\E^{|0,b]}}

    \def\mH{\mathcal{H}} \def\deF{de Finetti }
    \def\vt{\vartheta}
   
   \def\Le{Laplace exponent } \def\zrb{z^{b[}}
   \def\wrb{w^{b[}} 
       \def\mI{{\mathcal I}}

 \def\iffa{}  \def\eco{} \def\beco{} \def\np{}

\begin{document}

\selectlanguage{english}

\title[Spectrally negative Parisian L\'evy
processes,  and applications]{On fluctuation theory for
spectrally negative  L\'evy
processes with Parisian reflection below, and applications}

\author{Florin Avram}
\address{}
\curraddr{LMAP, Universit\'e de Pau,  France}
\email{florin.avram@orange.fr}
\thanks{.}

\author{Xiaowen Zhou}
\address{}
\curraddr{Concordia University, Montreal}
\email{xiaowen.zhou@concordia.ca}
\thanks{.}

\subjclass[2000]{Primary 60G51; Secondary 60K30, 60J75}
\date{28/11/2016}
\dedicatory{.}
\keywords{spectrally negative L\'evy process,  scale functions, capital injections, dividend optimization, valuation problem,    Parisian absorbtion and reflection}

\begin{abstract}
As well known, all functionals of a Markov process may be expressed  in terms of the generator operator, modulo some analytic work. In the case of spectrally negative Markov processes however, it is conjectured that
 everything can be expressed in a more direct way using the $W$ scale function  which intervenes in the two-sided first passage problem, modulo performing various integrals.
This conjecture arises from work on Levy processes \cite{AKP,Pispot,APP,Iva,IP, ivanovs2013potential,AIZ,APY}, where the $W$ scale function has explicit Laplace transform, and is therefore easily computable;
furthermore it was found in the papers above that
  a second  scale function  $Z$  introduced in \cite{AKP} (this is an exponential transform \eqref{Z} of $W$) greatly simplifies \fp laws, especially  for reflected processes.

   $Z$    is an  harmonic function of the L\'evy process (like $W$), corresponding to exterior boundary conditions $w(x)=e^{\th x}$ \eqref{hext}, and is also   a particular case of a "smooth \GS function" $\mS_w$. The concept of \GS function was introduced in \cite{gerber1998time}; we will use it however here in the more restricted sense of   \cite{APP15}, who define this
to be  a  "smooth"  harmonic function of the process, which fits the exterior boundary condition $w(x)$ and solves simultaneously the problems \eqref{serui}, \eqref{rserui}.

It has been conjectured  that  similar laws govern other classes of \sn processes,  but it is quite difficult to find assumptions which allow proving this  for general  classes of Markov processes.
However, we show below that in the particular case  of \sn \lev processes with Parisian absorption and reflection from below  \cite{AIZ,BPPR,APY}, this conjecture holds true, once the appropriate $W$ and $Z$ are identified (this observation seems new).

This paper gathers a collection of  first passage formulas for    spectrally negative Parisian L\'evy processes, expressed in terms of $W,Z$ and $\mS_w$, which may serve as an "instruction kit" for computing quantities of interest in applications, for example in risk theory and mathematical finance.
  To illustrate the  usefulness of our list,  we construct   a new index for the valuation  of
  financial companies modeled by spectrally negative L\'evy processes, based on a   Dickson-Waters modifications of the de Finetti optimal expected discounted  dividends objective. We offer as well an
  index for the valuation  of conglomerates of
  financial companies.


An implicit  question arising  is to investigate  analog results  for other classes of \sn Markovian processes.
\end{abstract}

\maketitle

 \sec{Introduction}

It is a great pleasure to have this opportunity to thank Nikolai Leonenko
for all the energy,  kindness and patience he put in our mathematical collaborations. He put things into my mind which will remain there, like for example the  Kolmogorov-Pearson-Wong diffusions, a  topic at the crossroads of probability and analysis.

As well known, the fluctuation/first passage theory of diffusions reduces
to the computation of  two non-negative monotone functions in the kernel of the associated Sturm-Liouville operator \cite{borodin2012handbook}, which are related to the \fp times $\t_b^{\pm}$ \eqref{tau}. A similar situation occurs for spectrally one sided L\'evy processes, when  one basic function,  the so-called $W$  scale function \cite{Ber}, suffices.

In the joint paper \cite{avram2009series}, we tried to put  these
two ``solvable models" together and study \fp problems for KPW  diffusions with one-sided jumps. The hope was to  unify the two theories, and I still hope to go back to this project one day.

The   paper below is a modest step in a different direction: it illustrates  the  fact that the fluctuation theory for
spectrally negative  L\'evy
processes with Parisian reflection below \cite{landriault2014insurance,
 AIZ,albrecher2015strikingly,PY,BPPR,APY} is formally identical with the classical one \cite{AKP,Pispot,APP,Iva,IP,
 ivanovs2013potential,AIpower,APP15}, once appropriate  scale functions are introduced.  Our hope is to stimulate further work on the intriguing question of whether this continues to be true for other classes of \sn Markovian processes. We also hope to provide an easily accessible summary of the general ideas, hoping to make more accessible this evolving literature
     with a huge potential for
 applications (mathematical finance, risk, inventory and queueing theory, reliability, etc).


Fluctuation theory (the study of maxima, minima and reflected processes) reduces in the  case of  spectrally one-sided L\'evy processes to the calculation of  scale functions. The name scale function in this context and the realization of the importance of the concept should be attributed to Bertoin \cite{Ber}. The idea was further developed in \cite{AKP,Pispot,APP,APP15} who illustrated that the answers to a wide variety of \fp problems may be
ergonomically expressed in terms of the so called $W,Z$ and $\mS_w$ scale functions, informally defined as the    harmonic functions of the process corresponding to exterior boundary conditions $0,e^{\th x}$ and $w(x)$, respectively. These functions are common to the free process $X$, to the process reflected below $X^{[0}$ or above $X^{b]},$ and to the doubly reflected process $X^{[0,b]}$, making thus natural to study these processes simultaneously.  Subsequently, \cite{Iva,IP,ivanovs2013potential} showed that the \lev formulas apply also in the much more general context of spectrally negative Markov additive processes,  once  appropriate
{\bf matrix scale}  functions are introduced.

Somewhat surprisingly,
we show here in section \ref{s:FPPar},  Proposition \ref{twostep}, that the classic formulas apply also to \sn \lev processes with Parisian absorbtion or  reflection below, once one identifies the appropriate scale functions $W,Z$
(this is true at least
 for the   \fp formulas we collected, and will probably be true for others).
In principle, similar formulas may hold for other classes of  \sn Markov processes, for example \lev processes with refraction \cite{KL,KPP,PY} or with a Parisian buffer \cite{APY}, or KPW diffusions with  jumps, and this topic deserves further investigation.

{\bf Contents}.
  We start with a  review in section \ref{LB} of
    some basics  \fp results for \sn L\'{e}vy processes, following \cite{AKP,IP,APP15}.

    We  sketch then in section \ref{s:Git}, as an  appetizer, an interesting financial application worth further study:    a  dynamic valuation index \eqref{e:git} for the constituents of a conglomerate of financial companies, which generalizes
the  \deF  optimal dividends objective \eqref{deFm}.
    Section \ref{s:MF} proceeds with a review of  some classic financial  {optimization problems}, solved using scale functions: the de Finetti and Shreve-Lehoczky-Gaver  optimization objectives, and the  American put option. These motivate their Parisian generalizations in the next
section \ref{s:FPPar}, the most important one, which gathers together eight
     recent \fp results for Parisian \sn L\'{e}vy processes \cite{AIrisk,AIZ,APY}. The goal is to provide a concise reference, and to
     illustrate the fact that  these results are {formally identical} with those for classic \sn L\'{e}vy processes, up to identifying the correct
     scale functions.

    In theorem \ref{t:eff}. section \ref{s:eff}  
    we   illustrate the usefulness of our list of formulas by calculating a static valuation index for financial companies modelled by Parisian \sn \lev processes, based on their "readiness to pay dividends".

 As a second application,  we provide in Section \ref{s:rei} a  dynamic valuation index for   "central branch networks" (conglomerates of financial companies with centralized decision taking).

 A typical example of  proof is included in Section \ref{s:pr}.

  \section{First passage theory  for \sn L\'evy processes, via the  scale functions\la{LB}}
   A  L\'evy process $X(t)_{t\geq 0}$ may be characterized   by its L\'evy-Khintchine/Laplace exponent,   defined by
\bea E_0 \big[e^{\th X(t)}\big]=e^{ t \k(\th)}.\eea

In applied probability we are  often confronted with the case of
\sn processes with negative jumps  only  (of course, \sp processes only require a change of sign), and with L\'evy-Khintchine decomposition of the form
 \beq \la{k} \k(\th) = \frac{\s^2}{2}\th^2 + p \th + \int_{\mbb R_+\backslash\{0\}}[e^{- \th y} - 1 + \th y ]\lm(d y), \;  \th \geq 0.\eeq
 Here the L\'evy measure (of $-X$)\fn[3]{Note that even though $X$ has only negative jumps, for
convenience we choose the L\'evy measure to have mass only on the positive instead of
the negative half line.} satisfies  $\lm(-\I,0)=0$ and
$ \int_{\mbb R_+\backslash\{0\}} (1\wedge y^2) \lm(dy) < \I$, and furthermore
   the drift (or profit rate) $p=E_0 [X(1)]=\k'(0_+)>-\I$ is finite (this requires
that the large negative jumps of the process have a finite mean, which is quite sensible, say in  modeling of catastrophes)\fn[4]{
 $X(t)$ is a Markovian process  with infinitesimal generator $\mG$, which acts
on $f\in C^2_c(\mbb R_+)$ as \cite[Thm. 31.5]{Sato}
\beq\label{eq:gamma11}
&&\mG h(x) = \frac{\s^2}{2}h''(x) + p h'(x) + \int_{\mbb R_+\backslash\{0\}}[h(x-y) - h(x) + yh'(x)]\lm(\td y);\no
\eeq
incidentally, this may be formally written as $\mG=\k(D)$, where $D$ denotes the differentiation operator.}.

A  further particular case to  bear in mind is when  the L\'evy measure has finite mass $\lm(0,\I)= \l<\I$. Writing $\lm(dz)= \l F(dz)$ allows
decomposing the resulting "Cram\'{e}r-Lundberg" process   as
  $$X (t)= x + c t -\sum_{i=1}^{N_\l} C_i, \; c=p+ \l \int_0^\I z F(dz).$$
  where $C_i, i=1,2,...$ are i.i.d. {nonnegative} jumps with distribution $F(dz)$,  arriving after
exponentially distributed times with mean $1/\l$.

 First passage theory  concerns the first passage times above and below, and the hitting time of a level $b$, defined by
\beq
 &\tau_b^+=\inf\{t\geq 0: X(t)>b\}, \qu \tau_b^-=\inf\{t\geq 0: X(t)<b\}, \qu \tau^{\{b\}}=\inf\{t\geq 0: X(t)=b\}
\la{tau} \eeq
 (with $\inf\emptyset=+\infty$). We write sometimes  $\t$ for the "ruin time" $\t_0^-$.

In applications, we are often interested  in versions  of  $X(t)$ which are constrained/regulated at first passage times below or/and above:
 \begin{eqnarray*}&& X^{[0}(t)= X(t) + R_*(t), \; \;  X^{b]}(t)= X(t)  -  R(t). \end{eqnarray*}
Here, \bea &&
 R_*(t)=-(\und X(t) \wedge 0), \; R(t)= R^b(t)=\le(\bar X(t)-b\ri)_+, \\&&  \und X(t)=\inf_{0 \leq s\leq t} X(t), \; \bar X(t):= \sup_{0 \leq s\leq t} X(t),\eea
  are  the minimal  Skorohod regulators constraining $X(t)$ to be nonnegative, and   to be smaller than $b$, respectively\fn[3]{Instead of $X^{b]}(t)$ one may study the "unused capacity" process $Y^b(t):=b -X^{b]}(t)=(\bar X(t)\vee b) -X(t)$}.
Also interesting  are processes    reflected at $0$ and refracted at the maximum (or  taxed) with coefficient $\d$, staring from $\t_b^+$ \cite{KL,KPP,PY}, \cite[(1)]{AIpower}:
  \begin{eqnarray*}&&X_\d(t)= X_\d^{[0,b[}(t)= X(t) + R_*(t)- \d R^b(t), \; \d \leq 1. \end{eqnarray*}
  The    regulators $R_*(t),R(t)$, defined by  having points of increase contained in $\{t\geq 0:X_\d(t)=0\}$ and $\{t\geq 0:X_\d(t)=\bar {X_\d}(t) \vee b\}$
   respectively, are  more complicated to describe explicitly    in this case; for a recursive construction, see
   \cite[Appendix]{AIpower}.

\beR
  The process $\d R(t)$ may be interpreted  as cumulative tax or dividends paid to some beneficiary.

 \eeR

\ssec{The  $W$  and $Z$ scale function for spectrally negative \lev processes}

Solving \fp problems involving
 Markovian  processes requires analytic work with their {\bf generator operator}. In the case of  L\'evy processes, the analytic work may be  replaced by the Wiener-Hopf factorization of the  { Laplace exponent} with  killing $\k(\th)-q$ (i.e. the identification and separation of its positive and negative roots).

  For  spectrally negative L\'evy process, the factorization involves only one \nne root
  \beq  \Phi_q:= \sup \{ \th \geq 0: \k(\th) - \q=0\},  \label{def_Fq} \eeq
  and everything reduces finally
 to the   determination of  one family of functions $W_q(x):\R -> [0, \infty), q \geq 0$,   defined on the positive half-line by the Laplace transform:
\beq \label{W}
&&\int_0^\infty  \mathrm{e}^{-\theta x} W_q(x) \diff x = \frac 1 {\k(\theta)-q}, \qu \for \theta > \Phi_q
\eeq
and taking the value zero on the negative half-line. Note that the singularity $\F_q$ plays a central role -- see for example \cite[Thm VII.1]{Ber}, stating that $\F_q$ is the \lev exponent of the subordinator $\t_x^+, x \geq 0$.

 The answers to many first passage problems may be ergonomically expressed in terms of this function, starting with the classic {\bf  gambler's winning problem}. For any $a < b$ and $x \in [a, b]$,
\beq &&
		\E_x \left( e^{-q \tau_b^+} 1_{\left\{ \tau_b^+ <  \tau_a^- \right\}}\right)  =\frac {W_q(x-a)}  {W_q(b-a)}=  \int_0^\I  e^{-q t} d P \left\{ \tau_b^+  <  \min[t,\tau_a^- ]\right\}.\fn[4]{Since this is the Laplace transform of the density of $\tau_b^+$, with absorbtion at
$\tau_a^-$, a  Laplace inversion will recover the corresponding density.}   		\label{twosided}
\eeq

\beR Note that establishing the equivalence between \eqref{twosided}
and \eqref{W}, "up to a constant", is not trivial. One solution when $q=0$,   via excursion theory, is provided in \cite[Thm VII.8]{Ber} by using the representation \beq W(x)=W(\I) e^{-\int_x^\I \fr{W'(h)}{W(h)} dh}= W(\I) e^{-\int_x^\I \nu(h) dh}, \la{exc}\eeq where $ \nu(h)=\fr{W'(h)}{W(h)}=n[\bar \e >h]$
is the characteristic measure of the Poisson process $n$ of excursion heights (informally, $\nu(h)$  is the rate of excursions starting at the moment $\t_h^+$ which are bigger than $h$, and thus cause ruin).
Later, \cite{nguyen2005some} used a Kennedy  type martingale, and \cite[(3)]{Pispot} constructed the scale function using potential theory
\beq \la{Wu} W_q(x)= \F'_q e^{\F_q x} -u_q(-x)=  u_{q}^{+}(-x) -u_q(-x)=e^{\F_q x}[u_q(0)-\fr{u_{q}^{+}(x) u_q(-x)}{u_q(0)}], \; x \geq 0, \eeq
where $\F_q$ is the inverse of the \lev  exponent \eqref{W} and $u_q$ is the potential density, which is exponential for nonnegative $x$, given by $u_{q}^+(x)=\F'_q e^{-\F_q x}, \; x \geq 0$.
\eeR
 {\bf A  second $Z_\q$ scale function}
  introduced  recently \cite{APP15,IP}\fn[6]{generalizing \eqref{oldZ} from \cite{AKP}}
 \begin{equation}Z_\q(x,\th)= e^{\th x} \le(1- \big(\k(\th)-\q\big) \int_0^x e^{ -\th y} W_\q(y) dy  \ri). \la{Z} 
 \end{equation}
   further simplifies  the solution  of many  first passage problems.

   For $\Re(\th)$ large enough {to ensure integrability},
it holds that \begin{eqnarray*}Z(x,\th)= \big(\k(\th)-\q\big)\int_0^\I e^{- \th y} W(x+y) dy. \end{eqnarray*}

Thus, the   $Z(x,\th)$ scale function is up to a constant an analytic extension of  the Laplace transform of the shifted scale function $W$ (also called
  Dickson-Hipp transform), and  the normalization ensures that  $Z(0,\th)=1$.
  $Z_\q(x,\th)$ is a  "smooth Gerber-Shiu function" with penalty $w(x)=e^{\th x}$ in the sense of \cite{APP15},
   i.e. the unique "smooth" solution of
\beq \la{hext}\bc (\mG-\q I) Z_\q(x,\th)=0, & x \geq 0\\ Z(x,\th)=e^{\th x},& x \leq 0\ec,\eeq
where $\mG$ is the Markovian generator of the process $X(t)$ -- see \cite[(1.12),(5.23), Sec 5]{APP15} and section \ref{nonh}. Indeed, consider the general solution $g(x,\th)=Z_\q(x,\th)+ k W_\q(x)$
and note that  continuity at $0$, i.e. $g(0,\th)=1 $ implies $ k=0$ when $W_\q(0) \neq 0$, and differentiability at $0$, i.e. $g'(0,\th)=\th $ implies $ k=0$ when $W_\q(0) = 0$ (since $W'_\q(0) > 0$).

Here are  some useful relatives of $Z_\q(x,\th)$ \cite{Kyp}:
\beq
&& \overline{W}_q(x):=  \int_0^{x} W_q(y) \diff y,   \no\\
&&Z_q(x) := Z_q(x, 0)=  1 + q \overline{W}_q(x)=\fr{\s^2}2 W'_q(y) + c W_q(y) -\int_0^x W_q(y) \ovl \lm(x-y) dy, \la{oldZ} \\ &&
\overline{Z}_q(x):= \int_0^{x} Z_q (z) \diff z = x + q \int_0^{x} \int_0^z W_q (w) \diff w \diff z
\eeq
{(the second definition of $Z_q(x)$ holds since $\mG \le(\overline{W}_q\ri)(x)=0$)}.
Note that  for $-\infty < x < 0$ we have
$    {Z_q(x, \F_q) = e^{\theta x}}, $ and  for $ x \leq 0$
$\overline{W}_q(x) = 0, \quad Z_q(x) = 1,  \quad \overline{Z}_q(x) = x$.

The $Z$ function  may be also characterized  by its  respective Laplace transform $\H Z_\q(s,\th)=
\fr 1 {\k(s)-\q}\fr{\k(s)-\k(\th)}{s-\th} $ (and $ \H Z_q(s)= \fr 1 {\k(s)-\q} \fr{\k(s) }{s}$).

Here are some examples of the utility of the $Z_\q(x,\th)$ function.
\beL \la{l:s}{\bf Severity of ruin for a process absorbed or reflected at $b>0$}.

 A)
 The joint Laplace transform of the first hitting time at $0$ and the undershoot
 is given by \cite[Prop 5.5]{APP15}, \cite[Cor 3]{IP}, \cite[(5)]{AIZ}
\beq \label{sevruin} &&
S_\q(x,\th):=\E_x\left(e^{-\q\tau_0^-}e^{\th X_{\tau_0^-}};\tau_0^-<\t_b^+\right) \no
\\&&=Z_\q(x,\th ) -   \fr {W_\q(x)}{ W_\q(b)}Z_\q(b,\th ), \; \theta \geq 0.
\eeq

   B) The joint Laplace transform of the first hitting time at $0$ and the undershoot  in the presence of reflection at a barrier $b \geq 0$ is
  \cite[Prop 5.5]{APP15}, \cite[Thm 6]{IP}
\beq\label{sevruinrefl} &&
S_\q^{b]}(x,\th)=\Eb_x\left(e^{-\q\tau_0^-}e^{\th X_{\tau_0^-}} \right)\no \\&&
=Z_\q(x,\th ) -   \fr {W_\q(x)}{ W'_\q(b)} Z'_\q(b,\th ), \;\theta \geq 0,
\eeq
where $\Eb$ denotes expectation for a process reflected from above at $b$, using Skorokhod reflection.

\eeL

\beR Note the similar structure of \eqref{sevruin} and \eqref{sevruinrefl};  formally, switching form absorbtion at $b$ to  the measure $\Eb$ involving reflection at $b$ only requires adding derivatives in the $b$ dependent coefficient.   This follows easily from the respective boundary conditions $S(b)=0,\le(S^{b]}\ri)'(b)=0$.

 \eeR

 \beR  By using
$
 \lim_{b\to\infty}\frac{Z_{\q}(b,\th)}{W_{\q}(b)}=
\frac{q-\k(\th)}{\Phi_{\q}-\th},$ we recover
$\E_x [e^{-\q \tau_0^- +\th X(\tau_0^-)} ]= Z_\q(x,\th) - W_\q(x)  \fr {\k(\th)-\q}{ \th-\Phi_\q} $ \;  \cite[(7)]{AIZ}. For $\th=0$ this is the famous ruin time transform
$\E_x [e^{-\q \tau_0^- } ]= Z_\q(x) - W_\q(x)  \fr {\q}{\Phi_\q} $. Note also the similar transform of the "recovery time" $\E_x [e^{-\q \t^{\{0\}} } ]=
\E_x [e^{-\q \tau_0^- +\Phi_\q X(\tau_0^-)}]=Z_\q(x,\Phi_\q) - W_\q(x)  \fr {\q}{\Phi_\q}$.
\eeR

 \begin{Rem}
The functions $Z_\q(x)$ and $Z_\q(x,\th)$ appeared first in \cite{AKP} and in  \cite[(5.23)]{APP15}, respectively. In the second paper (first submitted  in 2012)   $Z_\q(x,\th)$ was introduced as  a particular case of smooth Gerber-Shiu function associated to an {\bf exponential payoff} $e^{\th x}$,  as in Lemma \ref{l:s}, and as a generating function for the smooth Gerber-Shiu functions associated to  {\bf power payoffs}.

Subsequently, the ground-breaking papers \cite{IP,AIZ,AIpower} revealed several other \fp laws involving $Z_\q(x,\th)$.
 \end{Rem}

We turn now to a  result which we may be seen as  the fundamental  first passage law for reflected \sn \lev processes.
\beL {\bf The \LT \; of the discounted capital injections/bailouts until reaching an upper level}.
Let $X(t)$ denote a \sn L\'evy  process,
  let $R_*(t)=-(0\wedge \underline{X}(t))$ denote the regulator at $0$, let $X^{[0}(t)=X(t) +R_*(t) $  denote the process reflected at $0,$ and let $\E^{[0}_x$ denote expectation for the process   reflected at  $0$. The total capital injections into the  reflected  process, until the first up-crossing of a level $b$, satisfy
\cite[Thm 2]{IP}:
\begin{equation} \label{refbailout} 
L_*^{[0}(x,b):=\Ez_x [e^{-\q \t_b^+ - \th R_*(\t_b^+)}]=\bc Z_{\q}(x,\th) Z_{\q}(b,\th)^{-1}& \th <\I\\\E_x [e^{-\q \t_b^+} I_{ \t_b^+ < \t}]=
W_{\q}(x) W_{\q}(b)^{-1}&  \th =\I\ec,\end{equation}
where $W$ is the classic scale function \eqref{W}. \eeL
\beR \la{r:proof} This result  may be viewed as the fundamental law of \sn \lev processes, since it implies the smooth two-sided exit formula \eqref{twosided}. It is proved in  \cite[Thm 2]{IP} as a consequence of a more general result \cite[Thm 13]{IP}, but is essentially equivalent to  \eqref{sevruin}, by using
 \cite{IP}
  \beq \la{proof} \Ez_x [e^{-\q \t_b^+-\th R_*(\t_b^+)}]=
\E_x [e^{-\q \t +\th X(\tau_{0}^-)} ;  \t < \t_b^+ ] \; \Ez_0 [e^{-\q \t_b^+-\th R_*(\t_b^+)}]+ W_\q(x)W_\q(b)^{-1} \eeq
If the first term is known one gets an equation for the severity of ruin
\begin{eqnarray*}&& Z(x,\th) Z(b,\th)^{-1}= W_\q(x)W_\q(b)^{-1} + \E_x [e^{\th X(\tau_{0}^-)} ;  \tau_0^- <\t_b^+ ] Z(b,\th)^{-1},  \end{eqnarray*}
with the known solution
$
\E_x\left(e^{-\q\t + \th X(\t)};\t<\t_b^+\right)=Z_\q(x,\th) - W_\q(x) W_\q(b)^{-1} {Z_\q(b,\th)} $. And  if the severity of ruin is known, one may use \eqref{proof} with $x=0$ to solve for $\Ez_0 [e^{-\q \t_b^+-\th R_*(\t_b^+)}]$, provided that $W_\q(0) \neq 0$. When $W_\q(0) = 0$, one must start with a "perturbation approach", letting $x \to 0$  -- see for example Section \ref{s:pr}.
\eeR

Here is a generalization of Lemma \ref{l:s} B):
\beL   The {\bf  dividends- penalty law} for a  process reflected at $b$  is \cite[Thm 6]{IP}:
\begin{equation} \la{misL}
 S_{\vt}^b(x,\th):=\Eb_x\left[e^{ -q \t_0^-  + \th X(\t_0^-)-\vt \int_0^{\t_0^-} e^{-q t} dR(t)} \right]=
Z_\q(x,\th) - W_\q(x)    \fr{ Z_\q'(b,\th)+ \vartheta Z_\q(b,\th)  }{W_\q'(b)+ \vartheta W_\q(b)}.
\end{equation}
\eeL

\beR The function $S_\vt^b(x,\th)$ is $q$- harmonic, i.e. $(\mG- q)S_\vt^b(x,\th)=0$,  and satisfies $(S_\vt^b)'(b,\th)+ \vartheta S_\vt^b(b,\th)=0$.
Minimizing it when $\vt \neq 0$,
with respect to all possible dividend policies is an example of {\bf risk sensitive} dividends optimization.  It is expected that the last maximum of the barrier function
$G(b)= \fr{ Z_\q'(b,\th)+ \vartheta Z_\q(b,\th)  }{W_\q'(b)+ \vartheta W_\q(b)}$ will play a key role in the answer (note when $\th=\vt=0$, this is tantamount to maximizing the \deF barrier function $\fr{ W_\q(b)  }{W_\q'(b)}$).

\eeR

{\bf Extensions to other classes of \sn  Markov processes}.  Formulas with a similar structure in terms of  scale functions must hold for   other Markovian processes, and some recent papers  showed indeed that generalizations of the
$W,Z$ functions were still manageable computationally for other types of Markovian processes besides L\'evy: for example, see \cite{APP,LR,APP15} for reflected L\'evy processes, and \cite{KP,Iva} for spectrally negative Markov additive processes, and  \cite{AIZ,APY} for  \sn  L\'evy process with Parisian  reflection (bailouts). These results
increase considerably the arsenal of financial optimization tools.

\ssec{Nonhomogeneous problems and the Gerber-Shiu function $\mS_w$ \la{nonh}}
    When $e^{\th X({\tau_0^-})}$ is replaced  by an arbitrary penalty function $w( X({\tau_0^-})), w:(-\infty,0]\to\mbb R$  which is "admissible" (satisfies certain integrability condition) in \eqref{sevruin}, \eqref{sevruinrefl}, extensions of these formulas still hold if one   replaces $Z_\q(x,\th )$
by a
"smooth Gerber-Shiu function" $\mS_w$ \cite{APP15}.  More precisely,
given  $0<b<\infty$,  $x\in(0,b)$,
there exists a unique "smooth GS function" $\mS_w$ so that the following hold:
\beq&&
S_{w}(x )=\Eb_x\le[e^{-q
\t_a^-}w\le(X (\t_a^-)\ri)\mbf 1_{\{\t_a^-<\t_b^+\}}\ri]  = \mS_w(x ) - W_\q(x )
\frac{\mS_w(b ) }{W_\q(b )} , \la{serui} \\&&
S_w^{b]}(x )=\Eb_x\le[e^{-q
\t_a^-}w\le(X (\t_a^-))\ri)\ri]  = \mS_w(x )-
W_\q(x ) \frac{\mS_w'(b ) }{W_\q'(b )}.\la{rserui}
\eeq

   Stated informally, this amounts to the fact that both these  problems admit decompositions involving an identical "non homogeneous  solution" $\mS_w$.

The  "smoothness" required is:
 \begin{eqnarray}\label{dFw0}
\begin{cases}
\mS_w(0) = w(0),\\
\mS_w'(0_+) = w'(0_-), & \text{in the case $\sigma^2>0$ or $\lm([0,1])=\infty$}.
\end{cases}
\end{eqnarray}
Under these conditions, the function $\mS_w$ is unique and  may be represented as \cite[(5.13)]{APP15}:
\beq \la{GSdec} &&\mS_w(x)=   w(0) Z_\q(x) + \fr{\s^2}2 w'(0-)  W_\q(x) + \\&& \int_0^x W_\q(x-y) \int_{z=y}^{\I}  [w(0)-w(y-z)] \lm(dz)  dy \no\\&&=  w(0) \le(\fr{\s^2}2   W'_\q(x) + c W_\q(x)\ri)+ w'(0-) \fr{\s^2}2   W_\q(x) -   \int_0^{x} W_\q(x-y)  w^{(\lm)}(y) dy, \no \eeq
where $w^{(\lm)}(y)=\int_{z=y}^\I  [w(y-z)] \lm(dz)$ is the
 expected liquidation cost conditioned on a pre-ruin position of $y$ and on a ruin causing jump bigger than $y$.
 \beR The first two parts of \eqref{GSdec} may be viewed as boundary  fitting terms, and the last part as a "non-local Gerber-Shiu/nonhomogeneous component".  \eeR

\beL For  $w(x)=e^{\th x}$, the Gerber-Shiu function is $Z(x,\th)$ and the decomposition \eqref{GSdec} becomes:
\bea Z_\q(x,\th)= Z_\q(x) + \th \fr{\s^2}2 W_\q(x)+\int_0^x W_\q(y) \int_{x-y}^{\I}  [1-e^{\th(x-y-z)}] \lm(dz)  dy.\eea

This maybe easily checked by taking Laplace transforms, since
$$\H W(s)\fr{\k(s)-\k(\th)}{s-\th}=\H W(s)\Big(\fr{\k(s)}{s}+ \th \fr{\s^2}2+ \fr{\H \pi(s)-\H \pi(\th)}{s-\th}- \fr{\H \pi(s)-\H \pi(0)}{s}\Big)$$
\eeL

\sec{An appetizer: a dynamic index for valuation     of \sn \lev subsidiaries, built from the scale functions
\la{s:Git}}
\iffa
In this  section we introduce   a  valuation index
for subsidiaries. Note that in classic valuation in risk theory, by the \deF objective for example,  one fixes a prescribed liquidation point $o$  (usually $o=0$), in which case the index  can be expressed using the {\bf scale function methodology} started in \cite{AKP} and illustrated in this paper.
Here we propose however to further optimize the liquidation point, as is done for example in mathematical finance,
for valuing American put options.

Consider a subsidiary with {\bf liquidation value} $w(x)$ and   a policy $\pi=(R,R_*,\t)$ involving some dividend process $R$,   bailout process $R_*$, and
a  liquidation/reevaluation stopping time $\t$. If $d R_*(t) $ consists of a single payoff  at the liquidation time $\t$, then the {\bf modified De-Finetti objective} \eqref{e:robj} is
   \beq V^F(x)= {E_x \le[ \int_0^\t e^{- \q t} d R(t)  +   e^{- \q \t} w(X(\t))\ri]} \la{deFm}
  \eeq

  Subtracting now from the optimal value   obtained by  continuing until  some stopping time the immediate stopping value yields the following {\bf continuation index} alluded to in the title:
  \beq \mI(x)=\sup _\t {E_x \le[ \int_0^\t e^{- \q t} d R(t)
  +   e^{- \q \t} w(X(\t))\ri]- w(x)} \la{e:git} \eeq

\beR This is similar to the concept of Gittins index, in which we modify a value function by  subtracting  a constant subsidy $\mI=\mI(x)$
 for stopping, and choose this subsidy so that the {decisions} of whether to continue or stop yield equal payoffs.

 \eeR

  When the stopping time is prescribed by forced stopping, for example at $\t=\t_0^-$,
an explicit formula is available in terms of the   scale functions $W,Z$  and the \GS function $\mS_w$ \cite{LR,APP15}:
\beq &&V^F(x)=
 \mS_w(x)  +  W_\q(x)\fr {1-\mS_w'(b)}{ W'_\q(b)}, \la{deFmexp}\eeq

  With linear liquidation costs for example
$w(x)=\bc k x -K, & x<0\\x-K, &x \geq 0 \ec$,
the smooth \GS function is  \cite{LR,APP15}:
{ \beq  \mS_w(x)=k Z_{1,q}(x)-K Z_{q}(x), \; Z_{1,q}(x)=\bar Z_{q}(x)-p \bar W_{q}(x).
\eeq
}

  Optimizing \eqref{deFmexp} in $b$ is tantamount to optimizing   the  "barrier function"
\beq \la{G} G(b)=\fr {1-\mS_w'(b)}{ W'_\q(b)}. \eeq
Furthermore, the "barrier policy" at the last maximum $b^*$  is  often  optimal among all dividend policies, and, when barrier policies are not optimal,  an iterative procedure with starting point $b^*$ may be used to obtain the optimal "multi-bands policy"
\cite{schmidli2007stochastic,APP15}.

\beR
Several variations of the index \eqref{deFmexp} may be obtained  replacing absorption at $\t$ by Parisian absorbtion or reflection,  or by adding refraction
 or  other boundary mechanisms, which do not destroy the property of smooth passage upwards.

 Note that these problems do not require separate treatment:  the valuation index $\mI$ is again of the form  \eqref{deFmexp},  with  appropriate definitions of the  scale functions $W,Z$. Thus, finding two scale functions replaces the work necessary for many first passage problems   \cite{KL,AIZ,APY}.
\eeR

\beR The optimization  of the bailout point $o$ has been  less studied, and deserves further attention.\eeR

\sec{Some  applications of the $W$ and $Z$ functions to financial  optimization \la{s:MF}}
We turn now to reviewing applications of the $W$ and $Z$ functions to the optimization of several  financial   objectives involving  paying dividends and liquidation expenses, which seem  relevant for the problem of  evaluating the rentability/efficiency   of a subsidiary company.

{\bf The De Finetti objective with Dickson-Waters modification}
consists in  maximizing expected discounted dividends until the classic ruin time, amended  by  a  modification \cite{dickson2004some} penalizing the final liquidation:
   \beq\label{e:robj} &&V(x)=\sup_{\pi} E_x \le[ \int_0^{\tau_{0}^-} e^{- \q t} d R_{\pi}(t) +   e^{- \q \t_{0}^-} w(X(\tau_{0}^-))\ri]\no \\&& = V^F(x)+   S_{w}(x). \eeq
   Here $\pi$ represents an "admissible" dividend paying policy, $L^{\pi}(t)$ are the cumulative dividend payments, and $w(x)$ is a bail-out penalty function.
   The value function must satisfy, possibly in a viscosity sense,  the HJB equation \cite[(1.21)]{azcue2014stochastic}:
   \beq \mH ( V)(x):=\max[\mG_\q V(x),1-V'(x)]=0, x \geq 0, \qu V(x)=w(x), x < 0,\eeq
   where $\mG_\q V(x)$ denotes the discounted infinitesimal generator of the uncontrolled
surplus process, associated to the policy of continuing without paying dividends.
The second operator $1-V'_F(x)$ corresponds to the possibility of modifying the surplus by a lump payment.
The most important class of policies is that of constant barrier policies $\pi_b$, which modify  the surplus  only when  $X(t)>b$,
  by a lump payment bringing the surplus at $b$, and than  keep it there by Skorokhod reflection, until the next negative jump\fn[4]{in the absence of a Brownian component, this  amounts to paying all the income while at $b$}, until the next claim.

   Under  such a reflecting barrier strategy $\pi_b$, the dividend part of the De-Finetti objective
   has a simple expression in terms of the $W$ scale function:
   \begin{equation} V^F(x,b)=\Eazb_x\le[\int_{[0,\tau_{0}^-]}e^{-qs}d   R(s)\ri]=W_\q(  x)
W_{\q}^ {\prime}(b)^{-1}, \la{div} 
\end{equation}
where  $\Eazb$ denotes the law of the process  reflected from above at $b$,
and absorbed at $0$ and below.

The penalty part can be expressed in terms of $\mS_w(x)$;
finally, the modified de Finetti value function is:
\beq
V(x,b)=\bc \mS_w(x) +   W_\q(x)\fr {1-\mS_w'(b)}{ W'_\q(b)}& x \leq b\\x-b + V(b,b) &x \geq b  \ec  \la{VDeF}
\eeq
where $\mS_w$ is the smooth \GS associated to the penalty $w$ \cite{APP15}.

The "barrier function" \beq G_F(b):=\fr {1-\mS_w'(b)}{ W'_\q(b)}, b \geq 0 \la{GDeF}\eeq
  plays  
a central  role  in financial optimization.

 The most important cases of bail-out costs ${w}(x)$ are
\BEN \im exponential $w(x)=e^{\th x}$, when $\mS_w(x)=Z_\q(x,\th )$, and \im  linear  $w(x)= k x +K$. For $x < 0$, the constants $k >0$ and  $K\in \R$ may be  viewed as  proportional and  fixed bail-out costs, respectively.
The cases  $k \in (0, 1]$ and $k >1$   correspond  to  management being held responsible
for only part of the deficit at ruin, and to having to pay extra costs  at liquidation, respectively. When $K < 0$, late ruin is rewarded; when $ K> 0$ early liquidation is rewarded.

 In this case  $\mS_w(x)$ may be obtained by using
$Z_\q(x,\th )$ as generating function in $\th$, i.e. the coefficients of $K,k$ in $F(x)$ are found by differentiating  the $Z(x,\th)$ scale function  
$0$ and $1$ times respectively, with respect to $\th$. This yields
\beq \mS_w(x) =k \big(\ovl Z_\q(x) -p \ovl W_\q(x)\big) + K Z_\q(x ):=k Z_{1,\q}(x) + K Z_\q(x ). \la{Z1}\eeq
\EEN


\ssec{The Shreve-Lehoczky-Gaver infinite horizon objective, with linear penalties}. Suppose  a subsidiary must be bailed out each time its surplus is negative,
and  assume the penalty costs are linear  $w(x)=k x $. The optimization objective   of interest  combines   discounted dividends $R(t)$, and cumulative bailouts $R_*(t)$ 
  \beq\label{e:Sobj} V^S(x)=\sup_{b} \Ezb_x \le[ \int_0^\infty e^{- \q t} d R(t) + k \int_0^\infty e^{- \q t} d R_*(t)\ri]. 
  \eeq

  Since in a diffusion setting this objective has first been considered by Shreve, Lehoczky,  and Gaver (SLG) \cite{shreve1984optimal} -- see also  \cite{boguslavskaya2003optimization,lokka2008optimal} -- we will call it the SLG objective.

 The expected discounted dividends over an infinite horizon for the doubly reflected process,  with expectation denoted  $\Ezb$, are provided in   \cite[(4.3)]{APP}
\begin{eqnarray}\label{eq:DRL}
V^{S,D}(x,b)=\Ezb_x\le[\int_0^{\I} e^{-\q  t}d  R(t)\ri] &=&    Z_{\q}(x)  Z'_{\q}(b)^{-1}.
\end{eqnarray}

The capital injections part of the infinite horizon doubly reflected "SLG objective" is
 \beq  V^{S,w}(x,b)=\mS_w(x) -  \fr{Z_{\q}(x)}{ Z'_{\q}(b)} \mS_w'(b) \la{VSLG}\eeq

 \beR The Gerber-Shiu function $\mS_w(x)$ is common to  three
distinct nonhomogeneous problems involving a process $X$:
 \BEN \im severity of ruin with absorbtion at an upper barrier \eqref{serui} \im  severity of ruin with reflection  at an upper barrier \eqref{rserui} \im cumulative bailouts at the lower barrier with  reflection  at an upper barrier, for the doubly reflected process \eqref{VSLG}. \EEN
\eeR

\beR
 Note that in    \eqref{eq:DRL}, \eqref{VSLG}, just as in the relation $ \E_x e^{- \q \t_b^+}=Z_{\q}(x)  Z_{\q}(b)^{-1}$ \cite{AKP}, the second scale function $Z_{\q}(x)$ acts  for the process reflected at $0$ just as first scale function for the process absorbed at $0$. \eeR

In particular, with   linear costs  $w(x)=k x$, \eqref{VSLG} becomes:
  \begin{eqnarray*}
V^{S,w}(x)=-k \Ezb_x
\le[\int_0^{\I} e^{-\q  t}d  R(t)\ri] =
  k \le(Z_{1,\q}(x) -Z_{\q}(x)  Z'_{\q}(b)^{-1}  Z'_{1,\q}(b)\ri),
\end{eqnarray*}
 and  the optimal dividend distribution  is always of {\bf constant barrier} type \cite[(4.4)]{APP}.

 \ssec{Optimal dividend barrier strategies}

  {\bf The modified De Finetti dividend barrier function and the optimality of barrier strategies}. The last global maximum  $b^*$ achieving the $\max$ in
 $$\max_b V(x,b)=\mS_w(x) + W_{\q}(x) \fr{1- \mS_w'(b)}{ W'_{\q}(b)}$$
 plays a central role in  the optimal dividends distribution policy (even when this is not of single barrier type).  To determine $b^*$, it suffices to study the barrier  influence function
 $ G_F(b)=\fr{1- \mS_w'(b)}{ W'_{\q}(b)}$ \eqref{GDeF}.
 For example with linear costs    $w(x)=k x+K$,  this
  becomes
 \beq \la{GFlin} G_F(b)=\fr{1- k Z'_{1,\q}(b)-K Z'_{\q}(b)}{ W'_{\q}(b)},\eeq
 with $Z_{1,\q}(b)$ defined in \eqref{Z1};
  \cite[Lem. 4.1, Lem. 4.2]{LR} show that this barrier  function does attain a global maximum $b^*\in [0,\I)$, and that it is increasing-decreasing if $W'_{\q}(b)$ is log-convex.

\beXa With the SLG objective, the value function is
 $$\max_b V^{S}(x,b)=\mS_w(x) + Z_{\q}(x) \fr{1- \mS_w'(b)}{ Z'_{\q}(b)}$$
and the BF is $\fr{1- \mS_w'(b)}{ Z'_{\q}(b)}$.
With linear bailout costs    $w(x)=k x$,  this
  becomes
 \beq \la{GSlin} G_S(b)=\fr{1- k Z'_{1,\q}(b)}{ Z'_{\q}(b)}=\fr{1- k (Z_\q(b) -p  W_\q(b))}{\q W_{\q}(b)}.\eeq

After further removing a multiple of $ W_{\q}(b)$ from the numerator, we arrive at the equivalent optimization of
 $ \T G(b)=\fr{1- k Z_{\q}(b)}{\q W_{\q}(b)} $  \cite[(5.4)]{APP}.
\eeXa

\beR The  barrier functions and their largest maxima $b^*$ are easy to compute and central for solving numerically all barrier optimization problems, but determining whether the single barrier strategy at $b^*$  is optimal  is in many cases an open problem. However, the condition $G(0) \geq 0$   yields  simple criteria,
just  as in the case of the modified \deF objective; this motivated us to investigate also the optimality condition $G'(0) \leq 0$ -- see \cite{AM15,AM16} and section \ref{s:eff}. \eeR

\sec{A list of first passage laws for \lev processes with Poissonian (Parisian) detection of insolvency \la{s:FPPar}}

A useful type of models developed recently \cite{AIZ,albrecher2015strikingly,APY} assume  that insolvency is only {\bf  observed periodically},
  at an  increasing sequence of \emph{Poisson observation times} $\mT_r=\{t_i,i=1,2,...,$ the arrival times   of an independent Poisson process of rate $r$, with  $r > 0$ fixed\fn[4]{The concept of periodic observation may be   extended to the  Sparre Andersen  (non L\'evy) case,  using geometrically distributed  intervention times at the times of claims. This  deserves further investigation.}.

The analog concepts for first passage times are the stopping times
\begin{eqnarray}T_b^+= \inf\{t_i:\; X(t_i)  > b\}, \quad T_{a}^- = \inf\{t_i > 0:\; X(t_i) < {a}\} \label{def_tau_a_plus_minus}
 \end{eqnarray}
    We write sometimes  $T_0^-$ for  $T_0^-$.
 \iffa
  Under
 { Parisian observation times}, first passage is recorded
 only when the most recent excursion  below $a$/above $b$
 has exceeded an exponential rv $\kil_{r}$  of rate $r$.
\beR We will refer to stopping at $T_0^-$ as Parisian absorbtion.
 A {\bf \sn  \lev processes with Parisian reflection} below $0$ may be defined  by   pushing the process up to $0$ each time  it is below $0$ at an observation time $T_i$. In both cases, this will not be made explicit in the notation; classic and Parisian absorbtion and reflection will be denoted in the same way  (note that the first is a limit of the second).\eeR

\np

\beR    Parisian detection below $0$ is related to  the "time spent in the red"
$$ T_{<0}:= \int_0^\I I_{\{X(t) <0\}} d t, $$
 a fundamental risk measure studied by
 {\cite{picard1994some,loisel2005differentiation}.

  Indeed, the probability of Parisian ruin not being observed (and of recovering  without  bailout) when $p >0, q=0$ is   \cite[Cor 1,Thm 1]{LRZ11}, \cite[(11)]{AIZ}
   \beq  P_x[T_0^-=\I]=P_x[ T_{<0}< \kil(r)]
  = E_x \Big[e^{- r T_{<0}}\Big]= p \fr{\F_{r}} {{r}} Z(x,\F_r), \la{redext}\eeq}
 where $\F_{r}$  is the inverse of the \Le. This
 could be viewed as a {\bf state dependent} extension of the profit parameter $p$, to measure the profitability of risk processes.  {See \cite{li2014pre,li2015two}
 for further information about the relation between  Laplace transforms of occupation times
and the fluctuation theory of SNLP observed at independent Poisson arrival times}.

\eeR

\iffa

 The following proposition list  some basic \fp results for processes with Parisian detection of ruin, reflected or absorbed, following \cite{AIZ,BPPR,APY}. Note that these results coincide with the ones with classic, "hard" detection of ruin, and imply them when $r \to \I$. They also suggest that the known \fp results with hard ruin for SNMAPs \cite{KP,Iva,IP,AIrisk} might generalize to the Parisian case, provided that  properly defined scale matrix functions are introduced, and   written in  correct order.   To facilitate further work, we provide for each result below, besides the \lev Parisian reference, also the corresponding non-Parisian SNMAP reference from \cite{Iva,IP,AIrisk}, and also references from \cite{APP,AIZ,AIpower,albrecher2015strikingly} for  problems where  the SNMAP case is not yet resolved.
\beP \la{twostep} {\bf First passage results for processes with Parisian detection, followed by reflection or absorbtion}.
Let $X$ be a spectrally negative \lev process with Parisian detection below $0$, and fix $b>0$.   Assuming $x\in[0,b]$ and $\q, {r}>0,0\leq \th < \I$,  the  following hold:
\np
\BEN 
\im The {\bf   capital injections/bailouts law for a  Parisian reflected  process, until $\t_b^+$ \cite[Cor 3.1 ii)]{APY}}, \cite[Thm 2]{IP}. Let $X^{[0}(t)$ denote the SNMAP process reflected at $0,$
 let $R_*(t)=-(0\wedge \underline{X}(t))$ denote its regulator at $0$, so that $X^{[0}(t)=X(t) +R_*(t) $, and let $\Ez_x$ denote expectation for the process  with Parisian reflection at  $0$.  Then:
 {
\begin{equation} \la{Parisbailouts}  
B^b(x,\th):=\Ez_x [e^{-\q \t_b^+ - \th R_*(\t_b^+)}]=
\bc Z_{\q,r}(x,\th) Z_{\q,r}(b,\th)^{-1}& \th <\I\\E^{[0}[e^{-q\tau^+_b}; \tau^+_b<T_0^-]=
W_{\q,r}(x) W_{\q,r}(b)^{-1}&  \th =\I\ec,\end{equation}}
where
\begin{equation} \la{Z2}
 Z_{\q,{r}}(x,\th )=\fr{{r}}{\q+{r}- \k(\th)} Z_{\q}(x,\th)
 +\fr{\q-\k(\th)}{\q+{r}- \k(\th)} Z_{\q}(x,\Phi_{q+{r}} ),
 \end{equation}
 with  $\th=\Phi_{\q+r}$  interpreted in the limiting sense,
 and where
 ${W_{\q,{r}}(x)}:={Z_\q(x, \Phi_{\q+{r}})}$ \cite{AIrisk}, \cite[(12)]{AIZ}\fn[4]{The notation $W_{\q,{r}}(x):=Z_\q (x, \Phi_{\q+{r}})$  has been chosen to emphasize that this replaces, for processes with with Parisian ruin, the $W_\q$ scale function in  the classic "gambler's winning"  problem  $\E_x[e^{- \q \tau_b^+} ;\tau_b^+<\t ^-_0]=
\fr{W_{\q}(x)}{W_{\q}(b)}$.}.
When $r \to \I$, $Z_{\q,\I}(x,\th )=Z_{\q}(x,\th ), W_{\q,\I}(x )=W_{\q}(x)$ and \eqref{Parisbailouts}
reduces to classic results
 \cite{IP}.


 \im The  {\bf severity of Parisian ruin with absorbtion at $\t_b^+$},
   $S^b(x,\th)=\E_x\left[e^{\th X(T_0^-) }; 1_{T_0^- < \t_b^+ \wedge \kil_q}\right]$  is \cite[(15)]{AIZ} \cite[Cor 3]{IP}:
 {
\bea
&S^b(x,\th)=
Z_{\q,r}(x,\th)-W_{\q,r}(x)
{W_{\q,r}(b)}^{-1}{Z_{\q,r}(b,\th)}=Z_{\q,r}(x,\th)-L_{*,0}(x,b) {Z_{\q,r}(b,\th)}.
\eea
}

\im
Let $U_q^{|a,b|}(x,B)=\E_x \Big( \int_0^{ \tau^-_a \wedge \tau_{b}^+ } e^{-qt} 1_{\left\{ X(t) \in B  \right\}} \diff t\Big),$  denote the {\bf $q$-resolvent of a doubly absorbed \sn L\'{e}vy process with Parisian ruin \cite[Thm 2]{BPPR}}, for any Borel set $B\in[a,b]$. Then,
\begin{align} \label{resolvent_density} 
\begin{split}
	U_q^{|a,b|}(x,B) &= \int_{a}^{b} 1_{\{y \in B\}}\Big[ \frac {W_{\q,r}(x-a) W_{\q,r} (b-y)} {W_{\q,r}(b-a)} -W_{\q,r} (x-y) \Big] \diff y \quad a < y < b
	\end{split}
\end{align}
 (the analog classic result is Theorem 8.7 of \cite{K}).
\beR
It is natural to conjecture that the  resolvents for (partly) reflected processes
will also be of  the same form as the classic ones \cite[Thm. 1]{PDR},
\cite[Thm 2, Cor.2]{ivanovs2013potential}.
 \eeR

\im
 The {\bf  dividends- penalty law for a  process reflected at $b$, with Parisian ruin} \cite[Thm 6]{IP}, is:
\begin{equation} \la{mis1}
 S^b(x,\th,\vt):=\Eb_x\left[e^{  -\vt R(T_0^-)+ \th X(T_0^-)} ; T_0^- <\kil_q\right]=Z_{\q,r}(b,\th) - W_{\q,r}(b)   \fr{ Z_{\q,r}'(b,\th)+ \vartheta Z_{\q,r}(b,\th)  }{W_{\q,r}'(b)+ \vartheta W_{\q,r}(b)} \end{equation}
 \beq =
\big[Z_q(x,\th) -Z_q(x,\F_{q+r}) {H(b,\F_{q+r})}^{-1}{H(b,\th)}\big]{r}({r+q- \k(\th)})^{-1}, \la{mis}
\eeq
where $H(b,\th)= \vt Z_\q(b,\th)+ Z_\q'(b,\th)=(\th+ \vt) Z_\q(b,\th)-(\k(\th)-\q) W_\q(b)$\fn[5]{The structure of this formula reflects the fact that $\F_{q+r}$ is a removable singularity}.  The second, rather  complicated    formula,  is \cite[(23)]{AIZ}.

Note that when  $r\to\I$, \eqref{mis}   recovers the classic  \cite{AKP,nguyen2005some},\cite[Thm 6]{IP}, \cite[(25)]{AIZ}, by {using $
Z_q(b,\F_{q+r}) \to W_\q(b)$}:
\begin{eqnarray*} &&S^b_{\vt}(x,\th)=\Eb_x [e^{-\vartheta R(\t_0^-) + \th X(\t_0^-)} ; \t_0^- < \I]=\E _{Y(0)=b-x} [ e^{- \vt R(t_b) -\th (Y(t_b)-b)}] \no\\&&=Z_\q(x,\th) - W_\q(x) \le(W_\q'(b)+ \vartheta W_\q(b)\ri) ^{-1}  \le( Z_\q'(b,\th)+ \vartheta Z_\q(b,\th)  \ri)\fn[3]{This function satisfies $(S_\vt^b)'(b,\th)+ \vartheta S_\vt^b(b,\th)=0$.} \label{e:divtill0}\end{eqnarray*}
where
 $Y_x(t)= \overline{ X}(t) - X(t)$ is the  draw-down process/ reflection from the maximum.

At first sight,  \eqref{mis} and the classic version
\eqref{mis1}  look different; however, a little algebra will convince us that \eqref{mis} may also be written as \eqref{mis1}, with $Z_\q,W_{\q}$ updated to their Parisian versions $Z_{\q,{r}},W_{\q,{r}}$.

When $\vt=0$\fn[4]{When   $\vt=0=x$, and $X$ is L\'evy  with  bounded variation,
the {\bf  joint law of an excursion with reflection at an upper barrier, and of the final overshoot} is
\begin{eqnarray*}&&\Eb_0 [e^{ \th X(\t)} ; \tau_0^- < \kil_\q ]=1 + \fr 1 c  W_+'(b) ^{-1}  \le(W(b) \k(\th)- \th Z(b,\th)  \ri)      \end{eqnarray*}and $\Eb_0 [e^{ - \q \t} ; \tau_0^- < \I ]=1 - \T {\q}   \fr {W(b)} {W_+'(b)}, \;\;  \T \q =\fr {\q} c.        $
When $b \to \I,$  we recover the Laplace transform of an upward excursion $ \H \rho=\H \rho_\d=\E_0[ e^{-q \t}]=1-\fr{\T q}{\Phi_\q}= \fr {\l}c \H {\bar{F}}(\F( q)).$}
\eqref{e:divtill0} yields the {\bf severity of ruin for a regulated process} \cite{gerber2006note}:
\begin{equation} \Eb_x [e^{ \th X(\t)} ; \tau_0^- < \kil_\q ]=Z(x,\th) -  W(x) W'(b_+) ^{-1}  Z'(b,\th) =Z(x,\th) -  V^F(x,b)   Z'(b,\th). \label{regruinseverity}
\end{equation}

\beR
When $x=b$, we may factor the transform \eqref{mis} $\Eb_b\left[e^{ \th X(T_0^-) -\vt R(T_0^-)} ; T_0^- <\kil_q\right]$ as:
\begin{equation}
{\O}({\O+ \vt})^{-1}
\Bigg(Z_\q(b,\th)-{\O} ^{-1}\Big(\th Z_\q(b,\th) +\le(q-\k(\th)\ri) W_q(b)\Big)\Bigg) {r}({r+q- \k(\th)})^{-1}, \la{fac}\end{equation}
\bea \O=V^F(b,b)^{-1}=
{W_{q,r}'(b)}{W_{q,r}(b)}^{-1}=
{Z_q'(b,\F_{q+r})}{Z_q(b,\F_{q+r})}^{-1}=\F_{q+r}-r W_q(b) {Z_q(b,\F_{q+r})}^{-1}.\eea
Indeed, \bea &&Z_q(b,\th) -Z_q(b,\F_{q+r})\le( (\F_{q+r}+ \vt) Z_\q(b,\F_{q+r})-r W_\q(b)\ri)^{-1} H(b,\th)\\&&=Z_q(b,\th) -\le( \vt + \F_{q+r} - r W_\q(b) Z_q(b,\F_{q+r})^{-1}\ri)^{-1}H(b,\th)\\&&=Z_q(b,\th) -\le( \vt + \O\ri)^{-1} H(b,\th), \eea and
\eqref{fac} follows by simple algebra.
By \eqref{fac},  $ R(T_0^-)$ and $X(T_0^-)$ are independent when starting from $b$, and the former
has an exponential distribution with parameter $\O$ \cite[(23),(26)]{AIZ}.
\eeR

\im  The {\bf expected discounted dividends (upper regulation at  $b$)
 until $T_0^-$ \cite[(27)]{AIZ}}   are :
\begin{equation}
 V^{F}(x,b)=\Eb_x\le[\int_0^{T_0^-} e^{-\q  t}d  R(t)\ri] =   {W_{\q,{r}}(x)}{  W'_{\q,{r}}(b)}^{-1}.
\end{equation}

 \im  The {\bf expected discounted dividends
  with reflection at $0$ at Parisian times, until the total bail-outs  surpass an exponential variable $\kil_\x$  \cite[(15)]{AIpower}} are
 \begin{eqnarray}V^S(x,b,\th)= 
 \Ezb_x\le[\int_0^\I e^{-\q s} 1_{ [R_*(s) < \kil_\th]} d  R(s) \ri]={Z_{\q,r}(  x,\th)}
{Z_{\q,r}^ {\prime}(b,\th)^{-1}} \label{e:divskillbo}\end{eqnarray}

When $\th=0$, this becomes {\bf \cite[Cor 3.3]{APY}} \cite[(4.3)]{APP}:
\begin{eqnarray}\label{eq:DRLa2}
 V^S(x,b)=\Ezb_x\le[\int_0^{\I} e^{-\q  t}d  R(t)\ri] &=&    {Z_{\q,r}(x)}{  Z'_{\q,r}(b)}^{-1},
\end{eqnarray}
where $Z_{\q,{r}}(x)=Z_{\q,{r}}(x,0 ).$

\im
 The {\bf expected total discounted bailouts at Parisian times up to $\t_b^+$} are given for $0\leq x\leq b$ by
 \cite[Cor 3.2 ii)]{APY}:
\begin{align}\label{ParBail}
	V^F_*(x,b):=\Ez_x&\left[\int_0^{\t_b^+}e^{-qt}\diff R_*(t)\right]=Z_{\q,{r}}( x)  Z_{\q,{r}}( b)^{-1} \mS(b)-\mS(x).
\end{align}
where \beq \la{F} \mS(x)= \mS_{\q,{r}}(x)
=\frac{{r}}{\q+{r}}  \Big(\bar{Z}_q(x)+\frac{\k'(0_+)}{\q}
\Big).\eeq

\im  The {\bf total discounted bailouts at Parisian times over an infinite horizon}, with reflection at $b$ are \cite[Cor 3.4]{APY}:
\begin{eqnarray}\label{eq:localt0}
V_*^S(x,b)=\Ezb_x\le[\int_0^{\I} e^{-q t}\diff R_*(t)\ri] &=& Z_{\q,{r}}( x)  Z'_{\q,{r}}( b)^{-1} \mS'(b)-\mS(x).
\end{eqnarray}

\EEN

\eeP

\beR
Similar results hold
for  processes $X_\d^{b[}(t) $ with $\d$-refraction at a fixed point $b$ \cite{KL,KPP,Ren,PY}. The
scale functions are:
\begin{align}
\wrb_q(x)&=W_q(x)+\delta\int_b^x\mathbb{W}_q(x-y)W_q'(y)dy, \\
\zrb_q(x,\theta)&=Z_q(x,\theta)+
\delta\int_b^x\mathbb{W}_q(x-y)Z_q'(y,\theta)dy.
\end{align}

 For example, by \cite[Cor. 2]{KPP}, it holds that
\beq \la{nbr}
&&E_x \Big[e^{- {r} T_{<0}}\Big]=P_x[T_0^-=\I]= (p -\d)\fr{\F_{r}} {{r}- \d \F_{r}} \zrb(x,\F_r), \;  0\leq \d \leq p. \eeq
\eeR

\beR
Some of the results above have been extended   to processes $X_\d^{[0[}(t) $ with classic reflection at $0$ and refraction at the maximum \cite[(3)]{AIpower}. Thus, \eqref{Parisbailouts} holds with $Z_\q(x,\th)$ replaced by $Z_\q^{\fr 1{1-\d}}(x,\th)$ \cite[Thm 3.1]{AIpower}.  The proof uses the  probabilistic interpretation
$ \Ez_x [e^{-\q \t_b^+ - \th R_*(\t_b^+)}]=P[\t_b^+ < \kil_\q
\wedge K_{\th}],$
where $K_{\th}$ is the first time when the total bail-out exceeds
an independent exponential rv. $\kil_\th$.
Finally,  \cite[(22)]{AIZ}  extend this to the case when $\t_b^+$ is replaced by $T_b^+$.
\eeR

\beR The proof of the results above typically requires  in  the finite variation case only   applying the strong Markov property; however, in the infinite variation case,  the same problems require a perturbation approach-- see for example the proof of Proposition \ref{twostep} (A), in section \ref{s:pr}, or the use of the beautiful Ito excursion theory.
\eeR

\sec{Acceptance-rejection    of \lev subsidiary companies observed at Poissonian times, based on readiness to pay dividends\la{s:eff}}

Even in the one-dimensional case,
 the final choice of an acceptance-rejection principle  is not at all obvious.  A first intuition is that an  acceptable subsidiary must satisfy the classic positive profit condition \beq p:= E_0 [X(1)] >0\la{posload}\eeq
 or its extension involving  linear liquidation/bailout  costs  \cite{LR}.
However,  these equations only exploit the mean of the process involved and ignore its current state.
To remediate this deficiency, one could turn to model and state-dependent formulas like \eqref{redext}.  To be of practical use, an acceptance-rejection index should have  a complexity similar to that of the expressions above, and also intervene in some important optimization problem.

Note that the profitability/viability condition of \cite{LR} is equivalent to
$$G(0) \geq 0,$$
 where $G$ is the barrier influence function, and interesting variations may be obtained  by replacing  absorbtion at $0$ with reflection or Parisian reflection, which change the scale functions. The simplicity of all the resulting formulas
comes from the fact that the scale functions are only evaluated at $0$.
This suggested  an acceptance-rejection  criteria  introduced in  \cite{AM15,AM16}, based on the readiness of   subsidiaries to pay dividends at $b=0$.
   \beD \la{d:eff} A subsidiary will be called {\bf efficient} if the  barrier $b=0$  is  locally optimal  for paying dividends over some interval  $b \in [0,\e), \e>0,$ i.e. if it holds that $$G'(0) \leq 0.$$\eeD

   The motivation of this condition is that  companies satisfying it are functional even in the absence of cash reserves, and can contribute  cash-flows to the central branch without having to wait first until their reserves build out; efficiency is thus translated  as {\bf readiness to pay dividends}. This criterion turns out to be  a useful complement of the {\bf viability} concept $G(0) \geq 0$ (which at its turn generalizes the classic $p \geq 0$).
{\bf An additional bonus}  is that un-efficient  subsidiaries may be turned into efficient  ones by choosing an extra killing   $\q'_i$ to render  the barrier $b^*=0$
 locally  optimal; this means that the central branch will terminate subsidiaries deemed un-efficient by stopping  bailouts after times $\kil_{\q'_i}$ with exponential  law of  parameters $\q'_i$;  $(\q'_i)^{-1}$ will be referred to as   ``patience" parameters. The killing rate $\q'_i$  will be $0$ for subsidiaries deemed efficient.

 We illustrate now the application of the efficiency  concept for \sn L\'evy processes, under the SLG objective. 

  We assume  $\s=0,  \lm(0,\I)= \l <\I$,
  and      bail-outs at classic ruin times. This optimal dividend   problem is  fully analyzed in \cite[ Thm. 3]{APP}, and in particular
 \cite[Lem. 2]{APP}  show that  the optimal SLG constant barrier  is  $b^*=0$
iff
\beq   k \leq 1+ \fr {\q}{\l} \Eq \q=\l(k-1). \la{ksmall}\eeq

Note this is a simple application of the optimality of $0$ for the barrier function $\T G(b)=\fr{1- k Z_{\q}(b)}{\q W_{\q}(b)}$. Indeed,
\bea &&\q  \T G'(b)=-\fr{ k \q W^2_{\q}(b)+ W'_{\q}(b)(1- k Z_{\q}(b))}{ W^2_{\q}(b)}\\&&  \T G'(0) \leq 0 \Eq  k \q/ c^2+ (1- k ) (\q + \l)/c^2 \geq 0
   \Eq    \q + \l \geq k \l\eea

 \beR
  The  efficiency criterion \eqref{ksmall} does not take into account the law of $X$, beyond the total mass of its \lev measure.
However, it does have the interesting feature of making possible to turn  partially efficient subsidiaries 
into efficient ones, by introducing  extra killing $\q_i$. 
This follows from  the fact that
 the function $k(\q)$ which solves the equation $G'(0)=0$ is increasing. This   encouraging feature motivated us to remedy this by using  the SLG objective in more sophisticated environments including periodic observations, refraction, considered here, and also by using the De Finetti objective -- see \cite{AM16}. \eeR

The next result provides   a nontrivial  efficiency criteria under the  SLG infinite horizon cumulative  dividends-bailouts objective with  Parisian reflection
 \begin{Thm} \la{t:eff}

 a) The  SLG value function with  Parisian reflection and linear bailout costs $k x$ is:
\bea &&V_{SLG}(x)=Z_{\q,{r}}(x)  Z'_{\q,{r}}(b)^{-1}- k \le(Z_{\q,{r}}( x)  Z'_{\q,{r}}( b)^{-1} \mS'(b)-\mS(x)      \ri)
\\&& =k \mS(x)+Z_{\q,{r}}(x)  \fr{1- k \mS'(b)}{Z'_{\q,{r}}(b)}\eea

  b)
The barrier $b=0$ is a local maximum  iff the influence function $G(b):= \fr{1- k \mS'(b)}{Z'_{\q,{r}}(b)}$ satisfies
\beq && G'(0) \leq 0 \Eq k \Big ( \mS'(0)Z''_{\q,{r}}(0)-\mS''(0)Z'_{\q,{r}}(0)\Big)  \leq  Z''_{\q,{r}}(0) \no \\&&  \Eq k \fr{{r}}{\q+ {r}} \Big ( Z''_{\q,{r}}(0) -  \q  W_q(0_+) Z'_{\q,{r}}(0)\Big) \leq Z''_{\q,{r}}(0) \no \\&& \Eq
 k   \leq
(1+ \fr q r)\fr{\Phi_{q+r}-rW_q(0_+)}{\Phi_{q+r}-(r+\q)W_q(0_+)}.\la{effr}
\eeq

 In the  finite variation case\fn[4]{in the infinite variation case, the first equation still holds, but the efficiency index does not reflect the distribution, since $\Phi_{q+r}$ cancels} \eqref{effr}  holds   iff
{\small \bea &&k   \leq k(q,r):=
(1+ \fr q r)\fr{\Phi_{q+r}-r/c}{\Phi_{q+r}-(r+\q)/c} \la{effc0} \eea}
\end{Thm}
\beR

It may be checked  that $k(q,r)$ increases in $q$ from $k(0,r)=1$ to infinity and thus an inefficient  subsidiary with high  transaction cost $k > k(q,r)$ may be turned into efficient by increasing the killing $q$ sufficiently.
\iffa
More precisely, solving the inequality \eqref{effc0} yields $k < 1+ q/r$, or
{\small $$k \geq 1+ q/r, \; \; {q+r} < \k\le(\frac{(k-1) (q+r)}{c (k -1-q/r)}\ri) \Eq \F_{q+r} < \frac{(k-1) (q+r)}{c (k -1-q/r)},$$}
and inefficient subsidiaries may be made efficient by choosing an extra  killing rate $\q'_i$ obtained by letting $ \q_i:= \q + \q'_i$ solve the equality $\F_{ q+r} = \frac{(k-1)  ( q+r)}{c (k -1- q/r)}$. For example, for the \CL with   exponential claims and Laplace exponent
$\k(s)=s \le(s-\fr{\lambda}{\mu+s}\ri),$  this yields  a cubic equation in $ q$. The  solution
\bea  q=(k-1) \frac{\sqrt{(r+ \lambda +c \mu )^2+ 4 r \lambda( k-1-c \mu/r)}- \left(r+\lambda +c \mu \right)}{2 (k -1-c \mu/r )} \eea
is increasing in $k$, with a removable singularity  at $k=1+c \mu/r$.
\eco
The corresponding maximum of the barrier influence function is:
\bea G(0)=  \fr{1- k \fr{{r}}{\q+{r}}}{Z'_{\q,{r}}(0)}= \fr{1- k \fr{{r}}{\q+{r}}}{\frac{q}{q+r}\Phi_{q+r}}=
\fr{\q+{r}-  k {r}}{\q \Phi_{q+r}}\eea
\eeR

{\bf Proof of theorem \ref{t:eff}}.
{
 The SLG objective with Parisian bail-outs and its derivative is are:
\[
G(b)=\frac{1+k \mS'(b)}{Z_{q,r}'(b)}, \; G'(b)=\frac{kZ_{q,r}'(b) \mS''(b)-(1+k \mS'(b))Z_{q,r}''(b)}{(Z_{q,r}'(b))^2}.
\]

Efficiency is equivalent to
\beq \la{cor}
G'(0_+)\leq 0 \Lra k \le(Z_{q,r}'(0) \mS''(0)- \mS'(0))Z_{q,r}''(0)\ri)\leq Z_{q,r}''(0).
\eeq

Using
\bea &&Z'_{\q}(x, \th)=\th Z_{\q}(x, \th) +(\q-\k(\th)) W_{\q}(x),\\
&&Z''_{\q}(x, \th)=\th Z'_{\q}(x, \th) +(\q-\k(\th)) W'_{\q}(x)=\th^2 Z_{\q}(x, \th) +\th(\q-\k(\th)) W_{\q}(x)+ (\q-\k(\th)) W'_{\q}(x),\\
&&Z'_{\q,{r}}(x)=\fr{\q}{\q+{r}} \Big( \Phi_{\q+{r}} Z_{\q}(x, \Phi_{\q+{r}})-r W_\q(x)\Big) + \fr{{r}}{\q+{r}} \q W_\q(x)=\fr{\q}{\q+{r}}  \Phi_{\q+{r}} Z_{\q}(x, \Phi_{\q+{r}}), \\&& Z''_{\q,{r}}(x)=\fr{\q}{\q+{r}} \Phi_{\q+{r}}\Big( \Phi_{\q+{r}} Z_{\q}(x, \Phi_{\q+{r}})-  r  W_\q(x)\Big)
\eea
we find
\begin{align*}
Z_{q,r}'(b)&=\frac{q}{q+r}\Phi_{q+r}Z_q(b,\Phi_{q+r}),
Z_{q,r}''(b)=\frac{q}{q+r}\Phi_{q+r}\left(\Phi_{q+r}
Z_q(b,\Phi_{q+r})-rW_q(b)\right),\\
 \mS'(b)&=-\frac{r}{r+q}Z_q(b),
 \mS''(b)=-\frac{rq}{r+q}W_q(b),
\end{align*}
and finally
\begin{align*}
k\left(\Phi_{q+r}-(r+\q)W_q(0_+)\right)\leq
(1+ \fr q r)\le(\Phi_{q+r}-rW_q(0_+)\ri).
\end{align*}

The result follows by noting that the coefficient of $k$ is always positive \bsq

\beR In the finite variation case,  with $r \to \I,$ \eqref{cor} becomes
$k  \Big ( W'_{\q}(0) - \fr \q {c ^2}\Big) \leq W'_{\q}(0) \Eq k \leq 1+ \fr \q \l,$  which fits \cite[(5.6)]{APP}
(also,
$\lim_{r \to \I} \fr {  \Phi_{q+r}-r/c} {\Phi_{q+r}-(r+\q)/c }=1+\lim_{r \to \I} \fr {  \q/c} {\Phi_{q+r}-(r+\q)/c }=1 + \fr{\q}{ \l} $.
\eeR

\sec{Heuristic valuation  of CB networks,  using   claims line dividend policies \la{s:rei}}
We consider here one of the simplest risk networks, involving  a parent company/central branch, and several
  subsidiary branches \cite{AM15,AM16,ABPR}.
  \begin{Def} \la{d:CB } {\bf A
central branch (CB) risk  network  \label{ex:cb}} is formed from:
\BEN \im     Several \sn
subsidiaries  $X_i(t), i=1,\ldots,I$, which must be kept
  above  certain prescribed levels $o_i$ by bail-outs from a central branch (CB) $X_0(t)$, or be liquidated when they go below $o_i$.

  \im  The reserve of the CB  is a \sn process denoted by
 $X_0(t)$ in the absence of
  subsidiaries, and by $X(t)$ after subtracting the bailouts. The ruin time
$$
\t = \t_0^-=\inf\{t\ge 0: X(t) < 0\}
$$
causes the ruin of the whole network and leads to a severe penalty.

\im The CB  must also cover a certain proportion  $\bar \a_i=1 -\a_i$ of each  claim $C_{i,j}$ of subsidiary $i$, leaving the subsidiary to pay only  $\a_i C_{i,j}$, where $\a_i \in [0,1]$ are called  proportional reinsurance retention levels.

\EEN
\end{Def}

 A natural approach for {\bf evaluating financial companies}, going back to \deF \cite{de1957impostazione}  and Modigliani and Miller \cite{miller1961dividend} is  to   consider the optimal expected discounted cumulative dividends/optimal consumption until ruin\fn[4]{More generally, $\t$ could be replaced by other stopping times, like the drawdown (Azema-Yor) stopping time
$\t_\x:=\inf\{ t\geq 0: X(t) \leq \x \sup_{0\leq s \leq t} X(s) \},$
where $\x\in(0,1)$ is a fixed constant.}-- see \cite{leobacher2014bayesian} for further references on this venerable approach.

If the liquidation time $\t$ is also optimized
\beq \mI^F(\bff u):=\sup_{\pi=(R_0,R_1,...,R_I,\t)} E_{\bff u} \int_0^\t e^{- \q t} \le(\sum_{i=0}^I d R_i(t)\ri), \la{deFo}\eeq
the result $\mI^F(\bff u)$ is a Gittins type valuation index.

We will propose  here a heuristic multi-dimensional valuation index, based on a specific dividends policy, which,   remarkably, was found to be exact in \cite{AMP},  if the retention levels are small enough .

\np
We recall first from \cite{APP08,APP08a} that when $I=1$ and \bea 
  {c_0} \leq {c_1}\fr{\bar \a_1}{\a_1},\eea
  i.e. if the angle of the vector $\bff\a=(\a_1, \bar \a_1)$ with
the $u_1$ axis is bigger than that of $\bff c=(c_1,c_0)$, then
   the lower cone
$$\mathcal{C}:=\{0\leq {u_0} \leq {u_1}\fr{\bar \a_1}{\a_1}\}$$ contains $\bff c$ and is invariant  with respect to the stochastic flow,
\beco
i.e. that starting with initial capital
$(u_1,u_0)\in\mc C$, the process $(X_1,X_0)$ will stay there.
In particular, in the lower cone $\mathcal{C}$  ruin can only happen for the CB/reinsurer $X_0$,
and the \rp is a  classic one-dimensional ultimate ruin probability
$$\Rui(u_1,u_0)=\Rui(\a_1\fr{u_0} {\bar \a_1},u_0):=\Rui_0(u_0), \forall u_0$$
see for example \cite{Rolski,AA}.
Furthermore, the lower cone is also invariant with respect to the {\bf optimal} discounted dividends policy  \cite{AMP}.

Turning now to  several dimensions, it is easy to check that:
\beL The stochastic flow leaves invariant the cone
$$\mathcal{C}:=\{0 \leq u_0 \leq {u_i}\fr{\bar \a_i}{\a_i}, \; \bar \a_i= 1- \a_i, \; i=1,...,I\},$$ provided the "(extra) cheap reinsurance" condition
\beq c_0 \leq {c_i}\fr{\bar \a_i}{\a_i}, i=1,...,I\la{eqc} \eeq
is satisfied.

\eeL
The  boundary edge   \beq \la{equ} {u_1}\fr{1-\a_1}{\a_1}=...={u_i}\fr{1-\a_i}{\a_i}=
u_0, \; i =1,..., I,\eeq to be called "{\bf claims line}",
plays a prominent role in two recent
papers, \cite{bauerle2011optimal}\fn[4]{who computed an explicit value function  maximizing an expected exponential
utility at a fixed terminal time  for multi-dimensional reinsurance model under the "cheap reinsurance" assumption that the drifts point along the line
 $  {c_1}\fr{1-\a_1}{\a_1}=...={c_I}\fr{1-\a_I}{\a_I} $.} and \cite{AMP}, who solved
the optimal dividends problem  in the (extra) cheap reinsurance two-dimensional case
   $ {c_1}\fr{1-\a_1}{\a_1}> c_0.$
   The last paper showed that:
   \BEN \im Starting from the claims line, the optimal policy is to stay on this line by cashing the excess income of the subsidiary as dividends. \im  Starting from  points away from the claims line, in the  cheap reinsurance case, the optimal policy is to reach the claims line  by one lump sum payment.
   \im In the  extra cheap reinsurance case, the optimal policy is more complicated, when starting in a certain egg-shaped subset of the non-invariant cone (where   parts of the premia are cashed, following a "shortest path", in some sense).
\EEN
The first two findings prompt us  to introduce multi-dimensional "claims line" policies
for (extra) cheap reinsurance networks, under which the network follows this line in the absence of claims, by {\bf  subsidiaries cashing part of their premia as dividends}. Subsequently, whenever the CB or one subsidiary  drop below, {\bf all the other subsidiaries} reduce their reserves by {\bf lump sum dividend} taking, bringing back the process on the claims line.

\beR

These strategies may not be optimal; however, by postulating  that the   {subsidiary processes are just  linear functions  of the CB process,  they greatly simplify the problem}, and the value of the network expected dividends decomposes as a sum of  one-dimensional quantities-- see next Lemma.
\eeR

\beL For a general CB network, and a fixed  admissible  dividends process ${R_0}(t)$, the de Finetti value function for the equilibrium  policy associated to $\pi=({R_0},\t)$ is:
\begin{equation}
V^F_\pi(x)=E_x \Bigg[\int_0^{\t }e^{-qt} \Big[d {R_0}(t) + \T c dt -  \g d X_0(t) \no -\sum_{i=1}^I (\g \fr{\bar \a_i}{\a_i}-1) d X_i(t)\Big] \Bigg],
\end{equation}
where  $$\g= \sum_{i=1}^{I}\fr{\a_i}{\bar \a_i}, \; \T c= \g \sum_{i=1}^{I} c_i \fr{\bar \a_i}{\a_i}.$$
\eeL

Optimizing dividends reduces thus to a one-dimensional
problem.

\sec{Proof of Proposition \ref{twostep}.I)\la{s:pr}}

 By the Markov property,  we may decompose $g(x,b):=\E_x e^{-q\tau_b^+-\theta R_*(\tau_b^+)} , \theta>q+r,$ in three parts:
\begin{equation*}
\begin{split}
&g(x,b)
=\E_x[e^{-q\tau^+_b};\tau^+_b<\tau_0^-]
+\E_x\left[e^{-q\tau_0^-} \E_{X_{\tau_0^-}}[e^{-q\tau^+_0}; \tau^+_0<e_r];\tau_0^-<\tau^+_b \right]g(0,b)  \\ &+\E_x\left[e^{-q\tau_0^-}\E_{X_{\tau_0^-}}[e^{-qe_r+\theta X_{e_r}}; e_r<\tau^+_0];\tau_0^-<\tau^+_b  \right]g(0,b)\\
&=\frac{W_q(x)}{W_q(b)}+g(0,b)\left\{\E_x[e^{-q\tau^-_0+\Phi_{q+r}X_{\tau_0^-}}; \tau_0^-<\tau_b^+]+r\int_0^\infty e^{-\theta u} \right.\\
&\left.\left(W_{q+r}(u)\E_x[e^{-q\tau^-_0+\Phi_{q+r}X_{\tau_0^-}}; \tau_0^-<\tau_b^+] - \E_x[e^{-q\tau_0^-}W_{q+r}(X_{\tau_0^-}+u); \tau_0^-<\tau_b^+]   \right)du\right\}.
\end{split}
\end{equation*}
Noticing that
\begin{equation*}
\begin{split}
&\int_0^\infty  e^{-\theta u}E_x [e^{-q\tau_0^-}W_{q+r}(X_{\tau_0^-}+u) 1_{\tau^-_0<\tau^+_b}]du\\
&=E_x 1_{\tau^-_0<\tau^+_b}e^{-q\tau_0^-+\theta X_{\tau_0^-}}\int_0^\infty e^{-\theta(X_{\tau_0^-}+u)}W_{q+r}(X_{\tau_0^-}+u)du  \\
&=E_x 1_{\tau^-_0<\tau^+_b}e^{-q\tau_0^-+\theta X_{\tau_0^-}}\int_0^\infty e^{-\theta v}W_{q+r}(v)dv  =\frac{1}{\psi(\theta)-q-r}\left(Z_q(x,\theta)-
W_q(x)\frac{Z_q(b,\theta)}{W_q(b)}\right),
\end{split}
\end{equation*}
 we find
\begin{equation*}
\begin{split}
&g(x,b)=\left\{\frac{r}{\psi(\theta)-q-r}\left( Z_q(x,\Phi_{q+r})-Z_q(x,\theta)-W_q(x)\frac{Z_q(b,\Phi_{q+r})-Z_q(b,\theta)}{W_q(b)}\right)     +Z_q(x,\Phi_{q+r}) \right.\\
&\left. -W_q(x)\frac{Z_q(b,\Phi_{q+r})}{W_q(b)}\right\}g(0,b) +\frac{W_q(x)}{W_q(b)}=\le\{Z_{q,r}(x,\th)-W_q(x)
\frac{Z_{q,r}(b,\th)}{W_q(b)}\right\}g(0,b) +\frac{W_q(x)}{W_q(b)}.\\
\end{split}
\end{equation*}

Now in the finite variation case we may substitute $x=0$, and, using $W_q(0) >0,$ conclude that $g(0,b)=\frac 1{Z_{q,r}(b,\th)}$, which yields the result.

In the infinite variation case, we may use a  perturbation approach.
For $b>x>0$,  we have
\begin{equation}\label{g}
\begin{split}
g(0,b)
&=\E [e^{-q\tau^+_x}; \tau^+_x<e_r]g(x,b)+\E[e^{-qe_r+\theta X_{e_r}}; e_r<\tau^+_x, X_{e_r}<0]g(0,b) \\
&+\int_0^x \E[e^{-qe_r}; e_r<\tau^+_x, X_{e_r}\in dy]g(y,b)dy=e^{-\Phi_{q+r}x}g(x,b) +I_2(x)g(0,b)+I_3(x),
\end{split}
\end{equation}
\begin{equation*}
\begin{split}
I_2(x)&=r\int_{-\infty}^0\left(e^{-\Phi_{q+r}x}W_{q+r}(x-y)-W_{q+r}(-y)\right)e^{\theta y}dy\\
&=r\int_0^\infty e^{-\Phi_{q+r}x-\theta y}W_{q+r}(x+y) dy-\frac{r}{\psi(\theta)-q-r}\\
&=r\int_x^\infty e^{-\Phi_{q+r}x-\theta (z-x)}W_{q+r}(z) dz-\frac{r}{\psi(\theta)-q-r}\\
&=\frac{r}{\psi(\theta)-q-r}(e^{-\Phi_{q+r}x+\theta x}-1)- r\int_0^x e^{-\Phi_{q+r}x-\theta (z-x)}W_{q+r}(z) dz   \\
&=\frac{r}{\psi(\theta)-q-r}(e^{-\Phi_{q+r}x+\theta x}-1)+o(W_q(x)).
\end{split}
\end{equation*}

We can check that
\begin{equation*}
\begin{split}
&e^{-\Phi_{q+r}x}(Z_q(x,\Phi_{q+r})-Z_q(x,\theta))\\
&=e^{-\Phi_{q+r}x}\left[e^{\Phi_{q+r}x}(1-r\int_0^x e^{-\Phi_{q+r}y}W_r(y) dy)-e^{\theta x}(1-r\int_0^x e^{-\theta y}W_r(y) dy)  \right]   \\
&=1-e^{-\Phi_{q+r}x+\theta x}+o(W_q(x)),
\end{split}
\end{equation*}
\[Z_q(x,\Phi_{q+r})=e^{\Phi_{q+r}x}\left(1-q\int_0^x e^{-\Phi_{q+r}y}W_q(y)dy \right)  =e^{\Phi_{q+r}x}+o(W_q(x)), \text{   and}\]
\[ I_3(x)\leq \int_0^x E[e^{-qe_r}; e_r<\tau^+_x, X_{e_r}\in dy]dy=r\int_0^x e^{-\Phi_{q+r}x}W_{q+r}(x-y)dy=o(W_q(x)). \]

Solving now (\ref{g}) for $g(0,b)$ and letting $x\goto 0+$, we find again
\begin{equation*}
\begin{split}
g(0,b)&=\lim_{x\goto 0+}\frac{e^{-\Phi_{q+r}x} \frac{W_q(x)}{W_q(b)}}
{e^{-\Phi_{q+r}x} W_q(x)\frac{Z_q(b,\Phi_{q+r})}{W_q(b)}+re^{-\Phi_{q+r}x} W_q(x)\frac{Z_q(b,\Phi_{q+r})-Z_q(b,\theta)}{(\psi(\theta)-q-r)W_q(b)}+o(W_q(x))}  \\
&=\frac{\psi(\theta)-q-r}{(\psi(\theta)-q)Z_q(b,\Phi_{q+r})-rZ_q(b,\theta) }=\frac 1{Z_{q,r}(b,\th)}.
\end{split}
\end{equation*}

\end{document}